\documentclass[11pt]{amsart}
\usepackage{amsmath, amsthm, amscd, amsfonts}
\newtheorem{Lemma}{Lemma}

 \newfont{\Bbbb}{msbm10 scaled\magstephalf}
 \def\be{\begin{equation}}
 \def\ee{\end{equation}}

\def\bn{\begin{equation}}
\def\en{\end{equation}}
\def\br{\begin{center}}
\def\er{\end{center}}
\def\by{\begin{array}}
\def\ey{\end{array}}

\def\begy{\begin{eqnarray}}
\def\endy{\end{eqnarray}}
\def\bey*{\begin{eqnarray*}}
\def\eny*{\end{eqnarray*}}
\def\ber{\begin{tabular}}
\def\enr{\end{tabular}}

\def\bt{\begin{flushright}}
\def\et{\end{flushright}}

 \def\bea{\begin{equation}\begin{array}{ll}}
 \def\eea{\end{array}\end{equation}}

\begin{document}
\title[Difference of composition operators]{Compact differences of
composition operators from Bloch space to bounded holomorphic
function space in the Polydisc}
\author[Z.S. Fang and Z.H.Zhou]{Zhong-Shan Fang \and Ze-Hua Zhou$^*$ }
\address{\newline Department of Mathematics\newline
Tianjin Polytechnic University
\newline Tianjin 300160\newline P.R. China.}

\email{fangzhongshan@yahoo.com.cn}

\address{\newline Department of Mathematics\newline
Tianjin University
\newline Tianjin 300072\newline P.R. China.}
\email{zehuazhou2003@yahoo.com.cn}

\keywords{ Essential norm,Composition operator, Bloch space,
Polydiscs, Several complex variables}

\subjclass[2000]{Primary: 47B38; Secondary: 26A16, 32A16, 32A26,
32A30, 32A37, 32A38, 32H02, 47B33.}

\date{}
\thanks{\noindent $^*$Ze-Hua Zhou, corresponding author. Supported in part by the National Natural Science Foundation of
China (Grand Nos. 10671141, 10371091).}

\begin{abstract}
Let $\varphi$ and $\psi$ be holomorphic self-maps of the unit
polydisc $U^n$ in the $n$-dimensional complex space, and denote by
$C_{\varphi}$ and $C_{\psi}$ the induced composition operators. This
paper gives some simple estimates of the essential norm for the
difference of composition operators $C_{\varphi}-C_{\psi}$ from
Bloch space to bounded holomorphic function space in the unit
polydisc. Moreover the compactness of the difference is also
characterized.
\end{abstract}

\maketitle


\section{Introduction}

The algebra of all holomorphic functions with domain $\Omega$ will
be denoted by $H(\Omega),$ where $\Omega$ is a bounded domain in
$\mbox{\Bbbb C}^{n}.$ Let
$\varphi=(\varphi_1(z),\cdots,\varphi_n(z)) $ and
$\psi(z)=(\psi_1(z),\cdots,\psi_n(z))$ be holomorphic self-maps of
$\Omega.$ The composition operator $C_{\phi}$ induced by $\varphi$
is defined by
$$(C_{\phi}f)(z)=f(\phi(z)),$$ for $z$ in $\Omega$ and $f\in
H(\Omega)$.

We recall that the essential norm of a continuous linear operator
$T$ is the distance from $T$ to the compact operators, that is,
$$\|T\|_e=\inf\{\|T-K\|: K \mbox{ is compact}\}.$$ Notice that
$\|T\|_e=0$ if and only if $T$ is compact, so that estimates on
$\|T\|_e$ lead to conditions for $T$ to be compact.

During the past few decades much effort has been devoted to the
research of such operators on a variety of Banach spaces of
holomorphic functions with the goal of explaining the
operator-theoretic behavior of $C_{\varphi}$, such as compactness
and spectra, in terms of the function-theoretic properties of the
symbol $\varphi$.  We recommend the interested readers refer to the
books by J. H. Shapiro \cite{Sh1} and Cowen and MacCluer \cite
{Cow}, which are good sources for information on much of the
developments in the theory of composition operators up to the middle
of last decade.

In the past few years,  many authors are interested in studying the
mapping properties of the difference of two composition operators,
i.e., an operator of the form
$$T=C_{\varphi}-C_{\psi}.$$

The primary motivation for this has been the desire to understand
the topological structure of the whole set of composition operators
acting on a given functions. Most papers in this area have focused
on the classic reflexive spaces, however, some classical
non-reflexive spaces have also been discussed lately in the unit
disc in the complex plane. In \cite{Mac2}, MacCluer,Ohno and Zhao
characterized the compactness of the difference of composition
operator on $H^{\infty}$ spaces by Poincar\'{e} distance. Their work
was extended to the setting of weighted composition operators by
Hosokawa, Izuchi and Ohno \cite{HO1}. In \cite{Moo}, Moorhouse
characterized the compact difference of composition operators acting
on the standard weighted bergman spaces and necessary conditions on
a large scale of weighted Dirichlet spaces. Lately, Hosokawa and
Ohno \cite{HO1}and \cite {HO2} gave a characterization of compact
difference on Bloch space in the unit disc. In \cite{Toe} and
\cite{GorM}, Carl and Gorkin et al., independently extended the
results to $H^{\infty}(B_n)$ spaces, they described compact
difference by Carath\'{e}odory pseudo-distance on the ball, which is
the generalization of Poinar\'{e} distance on the disc.

The present paper continues this line of research, gives some simple
estimates of the essential norm for the difference of composition
operators induced by $\varphi$ and $\psi$ acting from Bloch space to
bounded function space in the unit polydiscs $U^n$, where
$\varphi(z)$ and $\psi(z)$ be holomorphic self-maps of the unit
polydisc in $n$-dimensional complex space. As its applications, a
characterization of compact difference is given .

\section{Notation and background}

Throughout this paper, let $D$ be the unit disc in the complex plane
$\mbox{\Bbbb C}$, $U^n$ the unit polydisc in $n$-dimensional complex
space $\mbox{\Bbbb C}^n$, and $|||z|||=\max\limits_j\{|z_j|\}$ stand
for the sup norm on $U^n$. For a holomorphic function in $U^n$,
define $\nabla f(z)=(\frac{\partial f}{\partial
z_1}(z),\cdots,(\frac{\partial f}{\partial z_n}(z)),$ $Rf(z)=<\nabla
f(z),\bar{z}>$.

For $z, w\in D$, the pseudo-hyperbolic distance between $z$ and $w$
is
 defined by $$\rho(z,w)=\left|\frac{z-w}{1-z\overline{w}}\right|.$$
 It is well known that  if $ f: D\rightarrow D$ is holomorphic, then
 $$\rho(f(z),f(w)\leq \rho(z,w)\hspace{2mm} z,w \in D.$$
The Bergman metric on the unit polydiscs is given by
 $$H_z(u,\overline{v})=\sum_{j=1}^n u_j\overline{v_j}/(1-|z_j|^2)^2.$$
The Kobayashi distance $k_{U^n}$ of $U^n$ is given by
$$k_{U^n}(z,w)=\frac{1}{2}\log
\frac{1+|||\phi_z(w)|||}{1-|||\phi_z(w)|||},$$ where $\phi_z:U^n
\rightarrow U^n$ is the automorphism of $U^n$ given by
$$\phi_z(w)=\left(\frac{w_1-z_1}{1-\overline{z_1}w_1},\cdots,\frac{w_n-z_n}{1-\overline{z_n}w_n}\right).$$

Let $H^\infty$ denote the space of bounded holomorphic functions $f$
on the unit polydiscs with the sup norm
$||f||_\infty=\sup\limits_{z\in U^n }|f(z)|$.

According to \cite{T1} and \cite{T2}, the Bloch space $\mathcal{B}$
in $U^n$ consists of those holomorphic functions such that
$$||f||_\mathcal{B}=\sup\limits_{z\in U^n}Q_f(z)<\infty$$
where $$Q_f(z)=\sup\left\{\frac{|\nabla f(z)\cdot
u|}{H^{\frac{1}{2}}_z(u,\overline{u})}: u\in \mbox{\Bbbb
C}^n-\{0\}\right\}.$$ It is well known that
$$||f\circ\phi||_\mathcal{B}=||f||_\mathcal{B}$$ for any automorphism $\phi$ of $U^n$, and
$\mathcal{B}$ is a Banach space under the norm
$$\|f\|_1=|f(0)|+\|f\|_\mathcal{B}.$$
If we put
$$G_f(z)=\sum_{j=1}^n(1-|z_j|^2)|\frac{\partial f}{\partial
z_j}(z)|$$ and $$||f||=|f(0)|+\sup\limits_{z\in U^n}G_f(z).$$ It
follows from $(3.8)$ and $(3.9)$ in \cite{T1} that
$$\frac{1}{n}G_f(z)\leq \max\limits_{1\leq j\leq
n}(1-|z_j|^2)|\frac{\partial f}{\partial z_j}(z)|\leq Q_f(z)\leq n
G_f(z),$$ this implies that $\frac{1}{n}\|f\|\leq \|f\|_1\leq
n\|f\|,$ so $\mathcal{B}$ is also a Banach space with the norm
$\|\cdot\|.$

\begin{Lemma} Assume $f\in \mathcal{B}$, then
$$|f(z)-f(w)|\leq n^2 \|f\|k_{U^n}(z,w)$$ for any $z,w \in U^n$.
\end{Lemma}
{\bf Proof.}
\begin{eqnarray*}
|f(z)-f(0)|&=&\left|\int_0^1\frac{Rf(tz)}{t}dt\right|=\left|\sum_{j=1}^n\int_0^1z_j
\frac{\partial f}{\partial \zeta_j}(tz)dt\right| \\
&\leq& \sum_{j=1}^n\int_0^1
\frac{|z_j|}{1-|tz_j|^2}\left|\frac{\partial
f}{\partial \zeta_j}(tz)\right|(1-|tz_j|^2)dt\\
&\leq&\|f\|_\mathcal{B}\sum_{j=1}^n\int_0^{|z_j|}\frac{1}{1-t^2}dt=\frac{1}{2}\|f\|_\mathcal{B}\sum_{j=1}^n
\log \frac{1+|z_j|}{1-|z_j|}\\
&\leq& n\|f\|_\mathcal{B} \frac{1}{2}\log
\frac{1+|||z|||}{1-|||z|||}.
\end{eqnarray*}
The last inequality follows by the fact the map $t\rightarrow
\log\frac{1+t}{1-t}$ is strictly increasing on $[0,1).$
 Setting $z=\phi_w(z)$, it follows that
\begin{eqnarray*}
|f(\phi_w(z))-f(0)|&\leq& n\|f\|_\mathcal{B} \frac{1}{2}\log
\frac{1+|||\phi_w(z)|||}{1-|||\phi_w(z)|||}\\
&=&  n||f\circ \phi_w||_\mathcal{B} \frac{1}{2}\log
\frac{1+|||\phi_w(z)|||}{1-|||\phi_w(z)|||}.
\end{eqnarray*}
That is,
$$|f\circ\phi_w(z))-f\circ\phi_w(w)|\leq n||f\circ \phi_w||_\mathcal{B} \frac{1}{2}\log
\frac{1+|||\phi_w(z)|||}{1-|||\phi_w(z)|||}.$$ Replacing $f\circ
\phi_w$ by $f\circ \phi_w \circ\phi_w^{-1}$,
$$|f(z)-f(w)|\leq n\|f\|_\mathcal{B}\frac{1}{2}\log
\frac{1+|||\phi_w(z)|||}{1-|||\phi_w(z)|||}\leq
n^2\|f\|k_{U^n}(z,w).$$ This completes the proof of the lemma.

\begin{Lemma} Suppose $f\in\mathcal{B}$, for fixed $0<\delta<1$, let $G=\{z\in U^n: |||z|||\leq \delta\}$. Then
$$\lim\limits_{r\rightarrow 1}\sup\limits_{||f||\leq 1}\sup\limits_{z\in G}|f(z)-f(rz)|=0.$$
\end{Lemma}
{\bf Proof}
\begin{eqnarray*}&&
\sup\limits_{z\in G}|f(z)-f(rz)|=\sup\limits_{z\in
G}\left|\sum\limits^{n}_{j=1}(f(rz_1,rz_2,\cdots,rz_{j-1},z_{j},\cdots,z_n)\right.\\
&&\hspace*{2cm}-\left.f(rz_1,rz_2,\cdots,rz_{j},z_{j+1},\cdots,z_n))\right|\\
&\leq& \sup\limits_{z\in
G}\sum\limits^{n}_{j=1}\left|\int^{1}_{r}z_{j}\frac{\partial
f}{\partial
z_{j}}(rz_1,rz_{j-1},tz_j,z_{j+1},\cdots,z^n)dt\right|\\
&\leq& (1-r)n\sup\limits_{z\in G}\left|\frac{\partial f}{\partial
z_{j}}(z)\right|\leq(1-r)n\sup\limits_{z\in G}\left|\frac{\partial
f}{\partial
z_{j}}(z)\right|(1-|z_j|^2)\frac{1}{1-|z_j|^2}\\
&\leq& (1-r)n\|f\|\sup\limits_{z\in G}\frac{1}{1-|||z|||^2}\\
&\leq& \frac{(1-r)n||f||}{1-\delta^2}.
\end{eqnarray*}
The lemma follows as $r\rightarrow 1$.

\section{Main theorem}

{\bf Theorem.}\hspace*{2mm} For $\delta>0$, write $F_\delta=\{z\in
U^n:\max(|||\varphi(z)|||,|||\psi(z)|||)\leq 1-\delta\}$. Suppose
$\varphi,\psi:U^n\rightarrow U^n$ and $C_\varphi-C_\psi: \mathcal{B}
\rightarrow H^{\infty}$ is bounded, then
$$\frac{1}{4}\lim\limits_{\delta\rightarrow 0}\sup\limits_{z\in E_\delta}|||\phi_{\varphi(z)}(\psi(z))|||
\leq ||C_\varphi-C_\psi||_e \leq 2n^2\lim\limits_{\delta\rightarrow 0}\sup\limits_{z\in E_\delta} k_{U^n}(\varphi(z),\psi(z))
$$ where $E_\delta=U^n-F_\delta.$

\begin{proof}
 We consider the upper estimate first. For fixed $0<r<1$, it easy
 to check both $C_{r\varphi}$ and $C_{r\psi}$ are compact
 operators. Therefore,
 $$||C_\varphi-C_\psi||_e\leq ||C_\varphi-C_\psi-C_{r\varphi}+C_{r\psi}||$$
Now for any $0<\delta<1$,
$$\hspace*{-4cm}||C_\varphi-C_\psi-C_{r\varphi}+C_{r\psi}||$$
\begin{eqnarray*}
&=&\sup\limits_{||f||\leq
1}||(C_\varphi-C_\psi-C_{r\varphi}+C_{r\psi})f||_\infty\\
&\leq&\sup\limits_{||f||\leq 1}\sup\limits_{z\in
F_\delta}|f(\varphi(z))-f(r\varphi(z))+f(r\psi(z))-f(\psi(z))| \\
&+&\sup\limits_{||f||\leq 1}\sup\limits_{z\in
E_\delta}|f(\varphi(z))-f(\psi(z))-f(r\varphi(z))+f(r\psi(z))|
\end{eqnarray*}
From Lemma 2, we can choose $r$ sufficiently close to $1$ such that
the first term of the right hand side is less than any given
$\epsilon$, and denote the second term by $I$. Using Lemma 1, it
follows that
\begin{eqnarray*} I&\leq & \sup\limits_{||f||\leq
1}\sup\limits_{z\in
E_\delta}(|f(\varphi(z))-f(\psi(z))|+|f(r\varphi(z))-f(r\psi(z))|)\\
&\leq&n^2\sup\limits_{||f||\leq 1}\sup\limits_{z\in
E_\delta}(k_{U^n}(\varphi(z),\psi(z))+k_{U^n}(r\varphi(z),r\psi(z)))\\
&\leq& 2n^2\sup\limits_{z\in E_\delta}k_{U^n}(\varphi(z),\psi(z))
\end{eqnarray*}
the last inequality follows by $k_{U^n}(r\varphi(z),r\psi(z)))\leq
k_{U^n}(\varphi(z),\psi(z)).$ First let $r\rightarrow 1$ and then
$\delta \rightarrow 0$, the upper estimate follows.

Now we turn to the lower estimate. Setting $$E_\delta^l=\{z\in U^n:
\max(|\varphi_l(z)|,|\psi_l(z)|)>1-\delta\}.$$ It is easy to see
that $E_\delta=\cup_{l=1}^nE_\delta^l$. For fixed $l$($1\leq l\leq
n),$ define
$$a_l=\lim\limits_{\delta\rightarrow 0}\sup\limits_{z\in E_\delta^l}
\frac{|\varphi_l(z)-\psi_l(z)|}{|1-\overline{\varphi_l(z)}\psi_l(z)|}.$$
If we put $\delta_m=\frac{1}{m}$, then $\delta_m \rightarrow 0$ as
$m\rightarrow \infty$.

If $||\varphi_l||_\infty=1$ or $ ||\psi_l||_\infty=1$, then for
enough large $m$ with $E^l_{\delta_m}\neq \emptyset$, so there
exists $z^m\in E_{\delta_m}^l$ such that $\lim\limits_{m\rightarrow
\infty}\frac{|\varphi_l(z^m)-\psi_l(z^m)|}{|1-\overline{\varphi_l(z^m)}\psi_l(z^m)|}=a_l$.
Since $z^m\in E_{\delta_m}^l$ implies that
$|\varphi_l(z^m)|>1-\delta_m$ or $|\psi_l(z^m)|>1-\delta_m$, without
loss of generality we assume $|\varphi_l(z^m)|\rightarrow 1$.
Setting
$$f_m(z)=\frac{1-|\varphi_l(z^m)|}{1-\overline{\varphi_l(z^m)}z_l}.$$
A little calculation shows that $\{f_m\}$ converges to zero
uniformly on compact subsets of $U^n$ as $m\rightarrow \infty$ and
$\|f_m\|\leq 2$ for any $m=1,2,\cdots$. So the compactness of $K$
implies that $||Kf_m||\rightarrow 0$ whenever $m\rightarrow \infty$,
it follows that
\begin{eqnarray*}
||C_\varphi-C_\psi-K||&\geq& \frac{1}{2}\limsup_{m\rightarrow \infty}||(C_\varphi-C_\psi-K)f_m||_\infty\\
&\geq& \frac{1}{2}\limsup_{m\rightarrow \infty}(||(C_\varphi-C_\psi)f_m||_\infty-||Kf_m||_\infty)\\
&=& \frac{1}{2}\limsup_{m\rightarrow \infty}||(C_\varphi-C_\psi)f_m||_\infty\\
&\geq&  \frac{1}{2}\limsup_{m\rightarrow \infty} \sup\limits_{z\in
U^n}|f_m(\varphi(z))-f_m(\psi(z))|\\
&\geq&  \frac{1}{2}\limsup_{m\rightarrow \infty} |f_m(\varphi(z^m))-f_m(\psi(z^m))|\\
&=&\frac{1}{2}\limsup_{m\rightarrow \infty}(1-|\varphi_l(z^m)|)
\left|\frac{1}{1-\overline{\varphi_l(z^m)}\varphi_l(z^m)}-\frac{1}{1-\overline{\varphi_l(z^m)}\psi_l(z^m)}\right|\\
&=&\frac{1}{2}\limsup_{m\rightarrow \infty}\frac{|\varphi_l(z^m)|}{1+|\varphi_l(z^m)|}
\left|\frac{\varphi_l(z^m)-\psi_l(z^m)}{1-\overline{\varphi_l(z^m)}\psi_l(z^m)}\right|\\
&=&\frac{1}{4}\limsup_{m\rightarrow
\infty}\left|\frac{\varphi_l(z^m)-\psi_l(z^m)}{1-\overline{\varphi_l(z^m)}\psi_l(z^m)}\right|\\
&=&\frac{1}{4}a_l=\frac{1}{4}\lim\limits_{\delta\rightarrow
0}\sup\limits_{z\in E_\delta^l}
\frac{|\varphi_l(z)-\psi_l(z)|}{|1-\overline{\varphi_l(z)}\psi_l(z)|}.
\end{eqnarray*}
If both $||\varphi_l||_\infty <1$ and $||\psi_l||_\infty<1$, in this
condition, when $\delta $ is small enough, $E_\delta^l$ is empty ,
without loss of generality, we may assume that
$$\lim\limits_{\delta \rightarrow 0}\sup\limits_{z\in
E_\delta^l}\left|\frac{\varphi_l(z)-\psi_l(z)}{1-\overline{\varphi_l(z)}\psi_l(z)}\right|=0.$$
Now for each  $l=1,2,\cdots, n$, we define
$$b_l=\lim\limits_{\delta\rightarrow 0}\sup\limits_{z\in
E_\delta}\left|\frac{\varphi_l(z)-\psi_l(z)}{1-\overline{\varphi_l(z)}\psi_l(z)}\right|$$
For any $\epsilon>0$, there exists a $\delta_0$ with $0<\delta_0<1$,
such that
$$\left|\frac{\varphi_l(z)-\psi_l(z)}{1-\overline{\varphi_l(z)}\psi_l(z)}\right|>b_l-\epsilon$$
whenever $z\in E_{\delta_0}$ and  $l=1,2,\cdots,n$. Since $z\in
E_{\delta_0}^l$ implies that $z\in E_{\delta_0}$, by the argument
above we have
\begin{eqnarray*}
||C_\varphi-C_\psi-K||&\geq& \frac{1}{4}\max\limits_{1\leq l\leq
n}\lim\limits_{\delta \rightarrow 0}\sup\limits_{z\in E_\delta^l}
\left|\frac{\varphi_l(z)-\psi_l(z)}{1-\overline{\varphi_l(z)}\psi_l(z)}\right|\\
&\geq&\frac{1}{4}\max\limits_{1\leq l\leq
n}(b_l-\epsilon)\\
 &=&\frac{1}{4}\lim\limits_{\delta \rightarrow
0}\sup\limits_{z\in
E_\delta}\max\limits_{1\leq l\leq n}
\left|\frac{\varphi_l(z)-\psi_l(z)}{1-\overline{\varphi_l(z)}\psi_l(z)}\right|-\frac{\epsilon}{4}\\
&=&\frac{1}{4}\lim\limits_{\delta\rightarrow 0}\sup\limits_{z\in
E_\delta}|||\phi_{\varphi(z)}(\psi(z))|||-\frac{\epsilon}{4}.
\end{eqnarray*}
Now the conclusion follows by letting $\epsilon \rightarrow 0$.
\end{proof}
{\bf Corollary} \hspace*{4mm}Suppose  $C_\varphi-C_\psi: \mathcal{B}
\rightarrow H^{\infty}$ is bounded, then $C_\varphi-C_\psi$ is
compact if and only if
$$\lim\limits_{\delta\rightarrow 0}\sup\limits_{z\in
E_\delta}|||\phi_{\varphi(z)}(\psi(z))|||=0.$$

\begin{proof} It follows from main theorem that the necessity is
obvious. By the strictly increment of $\log\frac{1+t}{1-t}$ on
$[0,1)$, $\lim\limits_{\delta\rightarrow 0}\sup\limits_{z\in
E_\delta}|||\phi_{\varphi(z)}(\psi(z))|||=0$ implies that
$\lim\limits_{\delta\rightarrow 0}\sup\limits_{z\in E_\delta}
k_{U^n}(\varphi(z),\psi(z))=0$, it follows from the main theorem
that $||C_\varphi-C_\psi||_e =0$, so $C_\varphi-C_\psi$ is compact,
the proof of this corollary is finished.\end{proof}

\end{document}